\documentclass[12pt]{amsart}
\usepackage{amssymb, epsfig,amsmath}
\usepackage[mathscr]{euscript}
\usepackage{amscd}
\usepackage{stmaryrd}

\newcommand{\bbp}{\mathbb{P}}
\newcommand{\bbq}{\mathbb{Q}}
\newcommand{\bbr}{\mathbb{R}}
\newcommand{\bbz}{\mathbb{Z}}
\newcommand{\bbf}{\mathbb{F}}

\DeclareMathOperator{\Gal}{Gal} \DeclareMathOperator{\End}{End}
\DeclareMathOperator{\GL}{\bf GL}

\newcommand{\gq}{\Gal(\overline{\bbq}/\bbq)}
\newcommand{\oqs}{\overline{\bbq}^{\sigma}}

\newtheorem{thm}{Theorem}[section]

\newtheorem{lem}[thm]{Lemma}
\newtheorem{prop}[thm]{Proposition}
\newtheorem{conj}[thm]{Conjecture}

\newtheorem{defnn}[thm]{Definition}

\newtheorem{remarkk}[thm]{Remark}

\newtheorem{examplee}[thm]{Example}

\title[Heegner points and Mordell-Weil groups  over large fields]
{Heegner points and Mordell-Weil groups of elliptic curves over
large fields}
\author{Bo-Hae Im}
\date{July 26, 2004}
\address{Department of Mathematics,
Indiana University, Bloomington, Indiana 47405}
\email{boim@indiana.edu} 

\begin{document}
\maketitle

\begin{abstract} Let  $E/\bbq$ be an elliptic curve
defined over $\bbq$ with conductor $N$ and $\gq$ the absolute
Galois group of an algebraic closure $\overline{\bbq}$ of $\bbq$.
We prove that for every $\sigma\in \gq$, the Mordell-Weil group
$E(\oqs)$ of $E$ over the fixed subfield of $\overline{\bbq}$
under $\sigma$ has infinite rank. Our approach uses the modularity
of $E/\bbq$ and a collection of algebraic points on $E$ -- the
so-called {\em Heegner points} -- arising from the theory of
complex multiplication. In particular, we show that for some
integer $r$, the rank of $E$ over all the ring class fields of
conductor of the form $rm$, and of the form $rp^n$  is unbounded,
as $m$ goes to infinity, as $m$ and as $n$ goes to infinity
respectively, where $m$ is a square-free integer and $p$ is a
prime such that $(m, rN)=1$ and $(p, rN)=1$.
\end{abstract}

\vspace{1 cm}
This paper is motivated by the following conjecture of M. Larsen
\cite{larsen}: \vspace{.5 cm}

 \noindent{\bf Conjecture.}\hspace{.1 cm} Let $K$ be a number field
and $E/K$ an elliptic curve over $K$. Then, for every $\sigma\in
\Gal(\overline{K}/K)$, the Mordell-Weil group
$E(\overline{K}^{\sigma})$ of $E$ over
 $\overline{K}^{\sigma}=\{x\in\overline{K}\mid \sigma(x)=x\}$
  has infinite rank.
\vspace{.5 cm}

 In \cite{im}, we have proved this conjecture in certain cases:

 for a number field $K$ and an elliptic curve $E/K$ over $K$,
 \begin{itemize}
 \item if 2-torsion points of $E/K$ are $K$-rational, or
 \item if  $E/K$ has a $K$-rational point $P$ such that $2P\neq O$ and $3P\neq O$,
   \end{itemize}

 then for every automorphism $\sigma\in
 \Gal(\overline{K}/K)$, the rank of  the Mordell-Weil group $E(\overline{K}^{\sigma})$
is infinite.

 In this paper, we prove that the conjecture is true for elliptic curves
 over $\bbq$ without any hypothesis on  rational points of
 $E/\bbq$, {\em i.e.} if $E/\bbq$ is an elliptic curve over $\bbq$, then, for every
 automorphism $\sigma\in
 \gq$, the Mordell-Weil group $E(\oqs)$ over the fixed subfield of $\overline{\bbq}$ under $\sigma$
  has infinite rank.

To prove the conjecture for a given $E/K$, ultimately one must
find an infinite supply of rational points of $E$ over finite
extensions of $K$ contained in $\overline{K}^{\sigma}$. In
\cite{im}, we constructed such points using Diophantine geometry,
essentially by searching for sufficiently rational subvarieties of
certain quotients of the $n$-fold product $E^n$ of $E$.

Here we use a completely different approach, coming from
arithmetic: taking advantage of the modularity of elliptic curves
over $\bbq$, we choose our rational points on $E$ to be algebraic
points over ring class fields -- the so-called {\em Heegner
points}. The main strategy is as follows:  by using the norm
compatibility properties of Heegner points and a generalized
dihedral group structure of the Galois groups of ring class fields
over $\bbq$,  we show that the rank of $E$ over the ring class
fields is unbounded as the ring class fields get larger. And we
also show that a given automorphism $\sigma\in \gq$ does not fix a
quadratic imaginary extension of $\bbq$ over which all primes
dividing the conductor of $E$ are split, then the rank of $E$ over
the fixed subfields of the ring class fields under $\sigma$ is
unbounded as the ring class fields get larger, which proves that
the rank of $E$ over $\oqs$ is infinite. On the other hand, if
$\sigma\in \gq$ fixes infinitely many quadratic imaginary
extensions, then we can construct infinitely many linearly
independent points defined over each of them by applying the
Hilbert irreducibility theorem to a Weierstrass equation of
$E/\bbq$ directly and this also proves infinite rank of $E$ over
$\oqs$.

%

\vspace{0.18 in}
\begin{center}
\sc{Acknowledgements}
\end{center}
\vspace{0.1 in}

I would like to thank my thesis advisor, Michael Larsen for
suggesting this problem and for valuable discussions and helpful
comments on this paper. Also, I wish to thank Henri Darmon for
suggesting this approach, and for his guidance and valuable advice
and  comments on an earlier manuscript of this paper.

\section{the main theorem}
In this section, we introduce the main theorem. First, we will
need the following Hilbert irreducibility and the denseness of
Hilbert sets in any open intervals of $\bbr$.

Let $f\in K(t_1,\ldots,t_m)[X_1,\ldots,X_n]$ be a polynomial with
coefficients in the quotient field $K(t_1,\ldots,t_m)$ of
$K[t_1,\ldots,t_m]$ which is irreducible over $K(t_1,\ldots,t_m)$.
We define $$H_K(f)=\{(a_1,\ldots,a_m)\in K^m\mid
f(a_1,\ldots,a_m,X_1,\ldots,X_n) \mbox{ is irreducible over }
K\}$$ to be the Hilbert set of $f$ over $K$.  If for every $m \geq
1$, any intersections of a finite number of  Hilbert sets with a
finite number of nonempty Zariski open subsets in $K^m$ are not
empty (in fact, they are infinite), a field $K$ is called \emph{a
Hilbertian field}.

\begin{lem}\label{lem:hilbert}
Let $L$ be a finite separable extension of a Hilbertian field $K$
and let $f$ be a polynomial in $
L(t_1,\ldots,t_m)[X_1,\ldots,X_n]$ which is irreducible over the
quotient field $L(t_1,\ldots,t_m)$. Then, there exists a
polynomial $p\in K[t_1,\ldots,t_m,X_1,\ldots,X_n]$ such that $p$
is irreducible over $K(t_1,\ldots,t_m)$ and $H_K(p) \subseteq
H_L(f)$.
\end{lem}

\begin{proof} For a given irreducible polynomial $f \in L(t_1,\ldots,t_m)[X_1,\ldots,X_n]$,
by (\cite{jar}, Ch.11, Lemma 11.6),
 there is an irreducible polynomial $q$  $\in
K(t_1,\ldots,t_m)[X_1,\ldots,X_n]$ such that $H_K(q)\subseteq
H_L(f)$. By (\cite{jar}, Ch.11, Lemma 11.1), there is an
irreducible polynomial $p \in K[t_1,\ldots,t_m,X_1,\ldots,X_n]$
which is irreducible over $K(t_1,\ldots,t_m)$ such that $H_K(p)
\subseteq H_K(q)$. Hence the Hilbert set $H_L(f)$ of $f$ over $L$
contains the Hilbert set $H_K(p)$ of $p$ over $K$.
\end{proof}

\begin{lem}
\label{lem:dense} Let $K$ be a number field and
$\tau_1,\ldots,\tau_m$ be a family of real embeddings of $K$. For
$i=1,2,\ldots,k$, let $f_i(x,y)\in K[x,y]$ be irreducible
polynomials over $K(x)$. Let $H_K(f_i)$ $= \{\alpha\in K\mid$ $
f_i(\alpha, y) \in K[y]$ is irreducible over $K\}$ be the Hilbert
set of $f_i$ over $K$. Then
$$\left(\bigcap\limits _{i=1}^k
H_K(f_i)\right)~\cap~\left(\bigcap\limits
_{j=1}^m\tau_j^{-1}(I)\right)\neq~~ \emptyset,$$ for any open
interval $I$ in $\bbr$.
\end{lem}

\begin{proof}
Since $K$ is a finite separable extension of  $\bbq$, by Lemma
\ref{lem:hilbert}, there exist irreducible polynomials
$F_i(x,y)\in \bbq[x,y]$ such that for each $i=1,2,\ldots,k$, the
Hilbert  set $H_{\bbq}(F_i)$ of $F_i$ over $\bbq$ is contained in
the Hilbert  set $H_{K}(f_i)$ of $f_i$ over $K$.

Let $I$ be an open interval in $\bbr$.  Since $\bigcap\limits
_{i=1}^k H_{\bbq}(F_i)$ is dense in $\bbq$ by (\cite{l83}, Chapter
9, Corollary 2.5), and $\bbq$ is dense in $\bbr$,
$\left(\bigcap\limits _{i=1}^k H_{\bbq}(F_i)\right)\cap I$ is not
empty. Hence there is $\beta \in \left(\bigcap\limits _{i=1}^k
H_{\bbq}(F_i)\right)\cap I$. Since $\bigcap\limits _{i=1}^k
H_{\bbq}(F_i)\subseteq \bigcap\limits _{i=1}^k H_K(f_i)$, we have
$\beta\in \bigcap\limits _{i=1}^k H_K(f_i)$. On the other hand,
for each real embedding $\tau_j$ of $K$, we have
$\tau_j|_{\bbq}=id_{\bbq}$. Hence for all $j=1,2,\ldots,m$,
$\tau_j(\beta)=\beta\in I$. Hence $\beta\in \bigcap\limits
_{j=1}^m\tau_j^{-1}(I)$. Therefore, $\beta\in \left(\bigcap\limits
_{i=1}^k H_K(f_i)\right)$ $\cap$ $\left(\bigcap\limits
_{j=1}^m\tau_j^{-1}(I)\right)$.
\end{proof}

Here is our main theorem.

\begin{thm}\label{thm:main} Let $E/\bbq$ be an elliptic curve
over $\bbq$. Then, for every automorphism $\sigma\in \gq$, the
rank of $E(\oqs)$ is infinite.
\end{thm}

\begin{proof} Let $N$ be the conductor of $E$ and let $y^2=x^3+ax+b$ be
a Weierstrass equation of $E/\bbq$. By the change of variables, we
may assume that $a$ and $b$ are integers. Let
$M=4p_1p_2\cdot\cdot\cdot p_k$, where $p_k$ are all distinct prime
factors of $N$. Consider the polynomial
$$f(x)=(1+Mx)^3+aM^4(1+Mx)+bM^6~~~\in\bbz[x].$$ Then, there
exists a real number $r$ such that for all $x<r$, the expression
$f(r)$ is strictly negative. Let $I=(-\infty, r)$ be the open
interval in $\bbr$ of all real numbers less than $r$. Since $\bbq$
is Hilbertian, by Lemma \ref{lem:dense}, there exists an integer
$m_1\in I$ such that $y^2-f(m_1)$ is irreducible over $\bbq$. Let
$K_{m_1}=\bbq(\sqrt{f(m_1)})$ be the quadratic imaginary extension
of $\bbq$. By Lemma \ref{lem:hilbert}, there is a polynomial $p$
over $\bbq$ such that $H_{\bbq}(p)\subseteq
H_{K_{m_1}}(y^2-f(x))$. Then, by Lemma \ref{lem:dense} again,
there exists an integer $m_2\in I\cap H_{\bbq}(p)$. So
$K_{m_2}=\bbq(\sqrt{f(m_2)})$ is a quadratic imaginary extension
of $\bbq$ and $K_{m_1}$ and $K_{m_2}$ are distinct, hence linearly
disjoint over $\bbq$. By repeating this procedure over the
composite field of quadratic imaginary extensions obtained from
the previous steps inductively, there is an infinite set $S$ of
integers  such that for all $m\in S$
\begin{enumerate}
\item $f(m)<0$, so that $K_m:=\bbq(\sqrt{f(m)})$ is a quadratic
imaginary extension of $\bbq$, \item the fields in the infinite
sequence $\{K_{m}\}_{m\in
S}$ are linearly disjoint over $\bbq$ 
and \item if $E/\bbq$ has CM, then, $K_m$ is different from
$F=\End(E)\otimes \bbq$.
\end{enumerate}

Note that for each $m\in S$ and for every prime $p_i$ dividing
$N$, $$f(m) \equiv \left\{\begin{array} {l@{}l} 1
\mbox{\hspace{.1in }} (\bmod~~p_i), \mbox{\hspace{.1in } if }
p_i\neq 2, \\ 1 \mbox{\hspace{.1in }} (\bmod~~8),
\mbox{\hspace{.1in } if } p_i= 2.
\end{array}\right.
$$

Hence, this implies that all primes dividing $N$ are split in
$K_m$. On the other hand, the discriminant of $K_m$ is $f(m)$ or
$4f(m)$ depending on whether $f(m)\equiv 1$ (mod 4) or not,
respectively. And in any case, the discriminant of $K_m$ is prime
to $N$.

Let $\sigma \in \gq$. We divide into two cases. For the first
case, suppose that for all $m\in S$, $\sigma|_{K_m}=id_{K_m}$.
Then, for each $m\in S$, consider the number $\displaystyle
\frac{1+Mm}{M^2}\in \bbq$. By plugging this number into the given
Weierstrass equation of $E/\bbq$, we get $$
y^2=\left(\displaystyle\frac{1+Mm}{M^2}\right)^3+a\left(\frac{1+Mm}{M^2}\right)+b
=\frac{f(m)}{M^6}.$$ Hence the point
$$P_m=\left(\displaystyle\frac{1+Mm}{M^2},
\frac{\sqrt{f(m)}}{M^3}\right)$$ is in $E(K_m)$ but it is not in
$E(\bbq)$.  And  moreover, since $K_m=K_m^{\sigma}$, $P_m$ is
fixed under $\sigma$.

So we get an infinite sequence $\{P_m\}_{m\in S}$ of points in
$E(\overline{\bbq}^{\sigma})$ such that each $P_m$ is defined over
the quadratic imaginary extension $K_m$ over $\bbq$. We may assume
that these points are not torsion points by (\cite{sil}, Lemma).
Now we show the points $P_m$ for $m\in S$ are linearly
independent. Suppose that for some integers  $a_i$ and $ m_i\in
S$,
$$a_1P_{m_1}+a_2P_{m_2}+\cdots+a_kP_{m_k}=O .$$ Since the fields in $\{K_m\}_{m\in S}$
  are linearly disjoint  over $\bbq$, for each $i$, there
  is an automorphism of $\overline{\bbq}$ which fixes all
  but one $K_{m_i}$ of $K_{m_1},\ldots,K_{m_k}$. Note that such an automorphism takes
  $P_{m_i}$ to its inverse, $-P_{m_i}$. Applying this
  automorphism, we get
  $$a_1P_{m_1}+\cdots +a_{i-1}P_{m_{i-1}}-a_{i}P_{m_{i}}+\cdots +a_kP_{m_k}=O.$$
By subtracting, we get $2a_iP_{m_i}=O$, which implies $a_i=0$. We
conclude that the $P_m$ for $m\in S$ are linearly independent in
$E(\oqs)\otimes \bbq$.
   Hence the rank of $E(\oqs)$ is infinite.

For the second case, suppose that there is an integer $m\in S$
such that $\sigma|_{K_m}\neq id _{K_m}$. Then, fix such a
quadratic imaginary extension $K_m$, and call it $K$, and let
$K_{ab}$ be the maximal abelian extension of $K$. Then, we
complete the proof of this case as a consequence of the following
stronger statement:

\begin{thm}\label{thm:ab} Under the assumption in the
second case above $($i.e. if $K$ is different from $\End(E)\otimes
\bbq$, all primes dividing $N$ are split in $K$, and
$\sigma|_K\neq id_K$$)$, the rank of the Mordell-Weil group of $E$
over the fixed subfield $(K_{ab})^{\sigma}$ of $K_{ab}$ under
$\sigma$ is infinite.
\end{thm}

The proof of Theorem \ref{thm:ab}, which lies deeper than the
methods used in the first case will be treated in the following
section and will be given explicitly in Proposition
\ref{prop:sigma}, as it requires modularity of $E/\bbq$ and the
theory of complex multiplication which give a non-torsion
algebraic point of $E$ under the given assumption in the second
case.
\end{proof}


%

\section{The rank of $E$ over ring class fields of imaginary quadratic fields}

The goal of this section is to prove Theorem \ref{thm:ab}
--- stated at the end of the previous section. By a theorem
of Wiles \cite{Wi} and Taylor-Wiles \cite{TW} (completed by a
later work of Breuil, Conrad, Diamond and Taylor \cite{BCDT}), the
elliptic curve $E/\bbq$ is known to be {\em modular}. Our strategy
is to construct algebraic points on $E((K_{ab})^{\sigma})$ using
{\em Heegner points} over the ring class fields arising from the
theory of complex multiplication.

We will need the following lemma later.

\begin{lem}\label{lem:prime1} Let $K$ be an
imaginary quadratic extension of $\bbq$ such that $K$ is different
from $F=\End(E)\otimes \bbq$. For a prime $\ell$, let $a_l$ be the
coefficient of the Hecke operator $T_{\ell}$ of the modular form
of $E$. Then, there is an integer $M$ such that for all $p\geq M$
$($which is inert in $F$, if $E$ has CM$)$, there is a prime $q$
such that
\begin{enumerate}
\item $q$ is inert in $K$, \item  $p$ does not divide $a_q$  and
\item $p$ divides $q+1$.
\end{enumerate}

\end{lem}

\begin{proof} Suppose $E/\bbq$ has no CM. Then, by (\cite{se}, \S4. Theorem 2),
 there is a large integer
$M$ such that for all primes $p\geq M$, the continuous Galois
representation $$\rho_p:\gq\rightarrow \GL_2(\bbf_p)$$ is
surjective. In particular, let $M$ be large enough such that every
$p\geq M$ is unramified in $K$, since $K$ is ramified at only
finitely many primes.

Since $Ker(\rho_p)$ is an open normal subgroup of $\gq$, there is
a finite extension $L$ over $\bbq$ such that
$Ker(\rho_p)=\Gal(\overline{\bbq}/L)$.

And since $$\Gal(L/\bbq)\cong\gq/Ker(\rho_p)\cong \GL_2(\bbf_p)$$
has a unique subgroup of index 2, the kernel of
$det^{\frac{p-1}{2}}$, the unique quadratic field $L'$ over $\bbq$
in $L$ is ramified only at $p$. On the other hand, since $p$ is
unramified in $K$, the fields $K$ and $L'$ are linearly disjoint.
Let $S$ be a finite set of primes such that the composite field
$KL'$ is unramified and $K$ is unramified outside $S$. Then, by
the Cebotarev density theorem, each Frobenius automorphism
Frob$_q$ comes up infinitely often in $\Gal(KL'/\bbq)$, for
$q\notin S$.

 Since
$K$ and $L'$ are linearly disjoint, we can choose the Frobenius
automorphism Frob$_q$ in $\Gal(K/\bbq)$ and $\Gal(L'/\bbq)$,
independently. Hence, there is a prime $q\neq 2$ such that $q$ is
inert in $K$ and
$$\rho_p(\mbox{Frob}_q)=
\begin{pmatrix} a& b\\0 &-a^{-1}\end{pmatrix}$$ where
 $a$ is in $\bbf_p$ such
that the order of $a$ is not $2$ in $\bbf_p$ and $b$ is some
element of $\bbf_p$. Then, we have
$$a_q=\mbox{Trace}(\rho_p(\mbox{Frob}_q))=a-a^{-1}\not\equiv
0~~~ (\bmod ~ p).$$  This implies that $p$ does not divide $a_q$,
which is  condition (2). Also, we have $$ q=
\mbox{det}(\rho_p(\mbox{Frob}_q))\equiv -1~~~ (\bmod~~ p).$$ This
implies that $p$ divides $q+1$ which is condition (3).
 Hence such a prime $q$ satisfies all conditions (1)
though (3).

Suppose $E/\bbq$ has CM. Let $R=\End(E)$ and $F=R\otimes \bbq$.
Then, $F$ is an imaginary quadratic extension of $\bbq$ such that
$R=\End_F(E)$.

By (\cite{se}, \S4. Corollary of Theorem 5), there is a large
integer $M$ such that for all primes $p\geq M$, the continuous
Galois representation $\rho_p$ factors through a surjective
homomorphism,
$$\sigma_p:\Gal(\overline{F}/F)\rightarrow R_p^*,$$ where
$R_p=R\otimes \bbz_p\cong \bbz_p^2$ as $\bbz_p$-modules. Under the
canonical embedding $R_p=\End_{R_p}(R_p)\subseteq
\End_{\bbz_p}(\bbz_p^2)=M_2(\bbz_p)$, the norm map and the trace
map from $R_p$ to $\bbz_p$ are the restrictions of the determinant
map and the trace map on $M_2(\bbz_p)$.

Let $\mathcal{O}$ be the maximal order in $F$. Then,
$\mathcal{O}/R$ is finite, so if $p\gg 0$, $$R/pR \cong
\mathcal{O}/p\mathcal{O}\cong \bbf_p\times\bbf_p \mbox{ or }
\bbf_{p^2},$$ depending on whether $p$ is split or inert in $F$.
In either case, the norm map $N$ is a non-degenerate quadratic
form on the underlying 2-dimensional vector space over $\bbf_p$.
For any non-zero $k$,  $N(\alpha)-kz^2$ is a non-degenerate
quadratic form on the 3-dimensional $\bbf_p$-vector space
$R/pR\times \bbf_p$. In particular, $N(\alpha)=-z^2$ defines a
conic curve. Since any conic has a rational point over a finite
field by the Chevalley-Warning theorem (\cite{se2}, Chapter I,
Theorem 3), the conic $N(\alpha)=-z^2$ over $\bbf_p$ is isomorphic
to $\bbp^1$, hence it has $p+1$ points. Since $N(\alpha)=0$ (when
$z=0$) has at most two non-trivial solutions on the line at
$\infty$, there are at least $p+1-3$ points of $N(\alpha)=-1$.
Therefore, for any prime $p\geq 5$, $N(\alpha)=-1$ has at least
three solutions over $\bbf_p$. On the other hand, the trace map
$T$ is a linear form on the underlying 2-dimensional vector space.
And the system of two equations, $T(\alpha)=0$ and $N(\alpha)=-1$
has at most two solutions. Since $N(\alpha)=-1$ has at least three
solutions for any prime $p\geq 5$,  we can find $\alpha\in R/pR$
such that
$$T(\alpha)\not\equiv 0 \mbox{~~and ~~} N(\alpha)\equiv -1
~(\bmod~p).$$

 Since the norm $N(\alpha)$ of $\alpha$ is congruent to
$-1$ (mod $p$), we can lift $\alpha$  to a $p$-adic unit
$\widetilde{\alpha}\in R_p^*$ by applying the Hensel's lemma to
the norm map.

Let $M\geq 5$ be large enough that for all $p\geq M$,
$R/pR\cong\mathcal{O}/p\mathcal{O}$ and let $p\geq M$ be inert in
$F$ and unramified in $K$.

 Let  $d_p$ denote the reduction of $det(\rho_p)$ mod $p$, and
let $\beta_p=d_p^{\frac{p-1}{2}}$. Then, $Ker(\beta_p)$ is of
index 2 in $\gq$ and therefore, corresponds to a quadratic
extension $L$ which is ramified only at $p$. Thus, the fields,
$F,K$ and $L$ are all linearly disjoint over $\bbq$. In
particular, $\sigma_p(\Gal(\overline{F}/FL))$ is of index 2 in
$R_p^*$, which has a unique open subgroup of index 2, since $p$ is
inert in $F$. As $FK\neq FL$, $K$ is not contained in the field
associated to $Ker(\sigma_p)$, which implies $K$ is linearly
disjoint from this field. By the Cebotarev density theorem, there
exists $q$ which is inert in $K$ and split in $F$ and such that
$$\sigma_p(\mbox{Frob}_q)=\widetilde{\alpha}\in R_p^*.$$
Then, we have
$$a_q=T(\sigma_p(\mbox{Frob}_q))=T(\widetilde{\alpha})
 \not\equiv 0 ~(\bmod ~p).$$ This implies
that $p$ does not divide $a_q$, which is  condition (2). Also, we
have $$ q= \mbox{det}(\sigma_p(\mbox{Frob}_q)) =
N(\widetilde{\alpha}) \equiv -1 ~(\bmod ~p).$$ This implies that
$p$ divides $q+1$ which is  condition (3). This completes the
proof.
\end{proof}

For a given elliptic curve $E/\bbq$, we fix an element
$\sigma\in\gq$ and a quadratic imaginary extension $K$ of $\bbq$
which is not fixed under $\sigma$. And we assume that all primes
dividing the conductor $N$ of $E/\bbq$ are split in $K$. Note that
this setting corresponds to the second case in the proof of
Theorem \ref{thm:main}.

 For each integer $n$ relatively
prime to $N$, let $H_n$ denote the ring class field of $K$ of
conductor $n$. By (\cite{dar}, Chapter 3, Theorem 3.6), there is
an algebraic point $P_n\in E(H_n)$ which is called a {\em Heegner
point of conductor $n$}. Let $HP(n)\subset E(H_n)$ denote the set
of all Heegner points of conductor $n$ in $E(H_n)$. Then, the set
$HP(n)$ satisfies the following properties.

First, we recall the norm-compatibility properties of Heegner
points.

\begin{prop}\label{prop:norm} Let $E/\bbq$ be a modular elliptic curve over
$\bbq$ and $N$ the conductor of $E/\bbq$. Let $n$ be an
integer and $\ell$ a prime number such that both $n$ and
$\ell$ are prime to $N$. Let $P_{n\ell}$ be any point in
$HP(n\ell)$ and $a_\ell$ the coefficient of the Hecke
operator $T_{\ell}$ of the modular form of $E$. Then, there
are  points $P_{n}\in HP(n)$ and $($when $\ell|n)$
$P_{n/\ell}\in HP(n/\ell)$ such that
$$\mbox{Trace}_{H_{nl}/H_n}(P_{nl})= \left\{\begin{array}
{r@{,\quad}l}  a_lP_n & \mbox{if } l\nmid n \mbox{ is inert in } K
\\ (a_l-\sigma_{\lambda}-\sigma_{\lambda}^{-1})P_n &\mbox{if }
l=\lambda\bar{\lambda}\nmid n \mbox{ is split in } K
\\ (a_l-\sigma_{\lambda})P_n & \mbox{if } l=\lambda^2 \mbox{ is ramified in } K
\\ a_lP_n-P_{n/l} & \mbox {if }
l|n. \end{array} \right.$$

\end{prop}
\begin{proof} See (\cite{dar}, Chapter 3, Proposition 3.10).
\end{proof}

\begin{lem}\label{lem:finite} Let $H_{\infty}$ be the union of all
the ring class fields of conductor prime to $N$. Then,  the set
$E(H_{\infty})_{tor}$ is finite.
\end{lem}

\begin{proof} See (\cite{dar}, Chapter 3, Lemma 3.14).
\end{proof}

The following lemma describes the structure of the Galois groups
of ring class field over an imaginary quadratic extension $K$ of
$\bbq$.

\begin{lem}\label{lem:degree} Let $N$ be the conductor of
$E/\bbq$. Let $H_n$ be the ring class field of conductor $n$ over
an imaginary quadratic extension $K$ over $\bbq$.

(a) If $p$ is a prime not dividing $c\cdot N\cdot [H_c:K]\cdot
disc(H_c)$, then for all integers $n\geq 1$,
$$ \Gal(H_{cp^{n}}/H_{cp})\cong \bbz/p^{n-1}\bbz,  \mbox{ and } \Gal(H_{cp^{n+1}}/H_{cp^n})\cong \bbz/p\bbz. $$

(b) If $k=\prod\limits_{j=1}^{m}p_j$ for distinct primes $p_j$
which are relatively prime to $N$ and inert in $K$, then for each
$j$,
$$ \Gal(H_k/H_{\frac{k}{p_j}}) \cong \bbz/(p_j+1)\bbz.$$
\end{lem}
\begin{proof} To prove the first part of (a), for $n\geq 1$ and $k\geq 1$,

$\Gal(H_{cp^n}/H_{cp}) \cong ker(\Gal(H_{cp^n}/K) \rightarrow
\Gal(H_{cp}/K))$

$\hspace{1.08in}\cong
ker((\mathcal{O}_K/cp^n\mathcal{O}_K)^*/(\bbz/cp^n\bbz)^*
\rightarrow (\mathcal{O}_K/cp\mathcal{O}_K)^*/(\bbz/cp\bbz)^*)$

$\hspace{1.08in}=
\displaystyle\frac{[(1+cp\mathcal{O}_K)/cp^n\mathcal{O}_K]^*\cdot(\bbz/cp^n\bbz)^*}{(\bbz/cp^n\bbz)^*}$

$\hspace{1.08in}=
\displaystyle\frac{[(1+cp\mathcal{O}_K)/cp^n\mathcal{O}_K]^*}{[(1+cp\bbz)/cp^n\bbz]^*}$

$\hspace{1.08in}=
\displaystyle\frac{[(1+p\mathcal{O}_K)/p^n\mathcal{O}_K]^*}{[(1+p\bbz)/p^n\bbz]^*},$
(since $p \nmid c\cdot[H_c:K]\cdot disc(H_c)$.)

$\hspace{1.08in}\cong
\displaystyle\frac{[(1+p\mathcal{O}_K)/p^n\mathcal{O}_K]^*}{(\bbz/p^{n-1}\bbz)},$

\hspace{1.08in} (since $[(1+p\bbz)/p^n\bbz]^*$ $ \cong
ker((\bbz/p^n\bbz)^* $ $ \rightarrow $ $ (\bbz/p\bbz)^*)) $ $\cong
\bbz/p^{n-1}\bbz.)$

 By the logarithmic function from
$1+p\widehat{\mathcal{O}}_{K,p} \rightarrow
p\widehat{\mathcal{O}}_{K,p}$ which maps
$1+p^n\widehat{\mathcal{O}}_{K,p} \mapsto
p^n\widehat{\mathcal{O}}_{K,p}$, where
$\widehat{\mathcal{O}}_{K,p}$ is the completion of
$\mathcal{O}_{K}$ at $p$,
$$[(1+p\widehat{\mathcal{O}}_{K,p})/p^{n-1}\widehat{\mathcal{O}}_{K,p}]^*
\cong
\widehat{\mathcal{O}}_{K,p}/p^{n-1}\widehat{\mathcal{O}}_{K,p}
\cong \mathcal{O}_{K}/p^{n-1}\mathcal{O}_{K}.$$ Hence, we have
$$\Gal(H_{cp^n}/H_{cp}) \cong (\mathcal{O}_{K}/p^{n-1}\mathcal{O}_{K})/(\bbz/p^{n-1}\bbz)\cong \bbz/p^{n-1}\bbz. $$

For the second part of (a),

$\Gal(H_{cp^{n+1}}/H_{cp^n})\cong ker( \Gal(H_{cp^{n+1}}/H_{cp})
\rightarrow \Gal(H_{cp^n}/H_{cp}))$

$\hspace{1.28in} \cong ker( (\bbz/p^n\bbz)^* \rightarrow
(\bbz/p^{n-1}\bbz)^*)$

$\hspace{1.28in} \cong \bbz/p\bbz.$

To prove (b), for each $i=1,...,m$,

$\Gal(H_k/H_{\frac{k}{p_i}}) \cong ker(\Gal(H_k/K) \rightarrow
\Gal(H_{\frac{k}{p_i}}/K))$

$\hspace{.96in} \cong
ker((\mathcal{O}_{K}/k\mathcal{O}_{K})^*/(\bbz/k\bbz)^*
\rightarrow
(\mathcal{O}_{K}/\frac{k}{p_i}\mathcal{O}_{K})^*/(\bbz/\frac{k}{p_i}\bbz)^*)$

$\hspace{.96in} \cong ker\left(\prod\limits_{j=1}^{m}
\Gal(H_k/H_{\frac{k}{p_j}}) \rightarrow \prod\limits_{j=1}^{i-1}
\Gal(H_k/H_{\frac{k}{p_j}}) \times \{1\}\times
\prod\limits_{j=i+1}^{m} \Gal(H_k/H_{\frac{k}{p_j}}) \right)$

$\hspace{.96in} \cong \Gal(H_k/H_{\frac{k}{p_i}} ) \cong
(\bbz/\lambda\bbz)^*/(\bbz/p_i\bbz)^*$, where $\lambda$ is a prime
factor of $p$.

 $\hspace{.96in} \cong \bbz/(p_i+1)\bbz$.
\end{proof}

{\em A Heegner system} attached to $(E,K)$ is a collection of
points $P_n\in E(H_n)$ indexed by integers $n$ prime to the
conductor $N$ of $E$, and satisfying the norm compatibility
properties given in (\cite{dar}, Chapter 3, Proposition 3.10
including Proposition \ref{prop:norm}) and the behavior under the
action of reflections described in (\cite{dar}, Chapter 3,
Proposition 3.11).

In our setting, since all primes dividing $N$ are split in $K$
(recall the construction of $K$ in the proof of Theorem
\ref{thm:main}), there is a Heegner system in which at least one
of the points $P_n$ for some $n$ is non-torsion by (\cite{dar},
Chapter 3, Theorem 3.13). We call such a Heegner system a {\em
non-trivial Heegner system}. We will need the following lemma.

\begin{lem}\label{lem:bdd}
Suppose there is  a nontrivial Heegner system attached to $(E,K)$.
Let  $n$ be a positive integer prime to the conductor of $E$ such
that 
there are non-torsion points of $E(H_{\infty})$ which are not in
$E(H_n)$. Then, there exists a positive integer $M$ such that for
all non-torsion points $Q\in E(H_{\infty})$ such that $Q\notin
E(H_n)$, $mQ\notin E(H_n)$ for all integers $m >M$.
\end{lem}

\begin{proof} Suppose not. Then, for all $M$, there exist an
integer $m>M$ and a non-torsion point $Q\in E(H_{\infty})-E(H_n)$
such that $Q, 2Q, \ldots, (m-1)Q \notin E(H_n)$ but $mQ\in
E(H_n)$.

We may assume that $m=p^k$ for some prime $p$ dividing $m$ and for
a positive integer $k$ by replacing $Q$ by $\frac{m}{p^k} Q$.
Then, either the exponent $k$ or the prime $p$ must go to infinity
as $M$ goes to infinity, since $m=p^k >M$.

Since $p^{k-1}<m$, $p^{k-1}Q\notin E(H_n)$. Hence there is an
automorphism $\tau\in$ Gal$({H_\infty}/H_n)$ such that
$\tau(p^{k-1}Q)\neq p^{k-1}Q$. Hence for the nontrivial point
$\tau(p^{k-1}Q)- p^{k-1}Q$,
$$ p(\tau(p^{k-1}Q)- p^{k-1}Q)=\tau(p^k Q)- p^kQ=p^k Q-p^k Q=O.$$
Since $p$ is a prime, this implies that the point $\tau(p^{k-1}Q)-
p^{k-1}Q$ has order exactly $p$, hence the nontrivial point
$\tau(Q)- Q$ has order exactly $m=p^k$. So we have shown that
there are torsion points of order exactly $m$ for all integers $m
> M$. Therefore, as $M$ goes to infinity, the order of torsion points in
$E(H_{\infty})$ is unbounded, which is a contradiction to the
finiteness of the set of torsion points of $E(H_{\infty})$ in
Lemma \ref{lem:finite}. This completes the proof.
\end{proof}

Now,  by using Lemma \ref{lem:prime1} and the norm-compatibility
properties in Proposition \ref{prop:norm}, we prove that the
subgroup generated by all of Heegner points of the given
nontrivial Heegner system is not finitely generated as a subgroup
of the elliptic curve over the union of all the ring class fields
of conductor of the form $rm$ for some $r$, where $m$ is a
square-free integer relatively prime to $rN$. Then, by using Lemma
\ref{lem:bdd}, we show the following unboundedness of the rank of
Mordell-Weil groups over all of those ring class fields.

\begin{prop}\label{prop:inf} Let $N$ be the conductor of the given
elliptic curve $E/\bbq$. If there is a nontrivial Heegner system
attached to $(E,K)$, where $K$ is a quadratic imaginary  extension
of $\bbq$ which is different from the field $\End(E)\otimes \bbq$
such that there is a non-torsion Heegner point $P_{r}\in E(H_{r})$
over some ring class field $H_{r}$ of conductor $r$, then the rank
of $E(H_{rm})$ is unbounded, as a square-free integer $m$ such
that $(m, rN)=1$ goes to infinity.
\end{prop}

\begin{proof} First, we show that  the group generated by all of Heegner
points of conductor $rm$ for all square-free integers $m$ such
that $(m,rN)=1$ in the given nontrivial Heegner system is not
finitely generated. Suppose this group is finitely generated.
Then, since there are only finitely many torsion (Heegner) points
over all the ring class fields in the system by Lemma
\ref{lem:finite}, for some integer $n$ which is a square-free
multiple of $r$ and $(n/r, rN)=1$, there is a fixed ring class
field $H_{n}$ over which all Heegner points of conductor $rm$ for
all square-free integers $m$ are defined. And we may assume that a
Heegner point $P_{n}$ of level $n$ is of infinite order by the
assumption that there exists a non-torsion Heegner point $P_{r}\in
E(H_r)$.

 Only finitely
many primes divide either $P_{n}$ in $E(H_{n})$ or points of
$E(H_{\infty})_{tor}$, because $E(H_{n})$ is finitely generated by
the Mordell-Weil Theorem and there are only finitely many torsion
points over all the ring class fields in the system by Lemma
\ref{lem:finite}. Let $S$ be the finite set of primes which divide
either $P_{n}$ in $E(H_{n})$ or the order of any point of the
finite set $E(H_{\infty})_{tor}$. Then, we can choose a large odd
prime $p\notin S$ such that $p\geq n+2$ and if $E$ has no CM, then
$p$ is unramified in $K$ and if $E$ has CM, then $p$ is inert in
the imaginary quadratic extension, $\End(E)\otimes \bbq$ of
$\bbq$.

Then by Lemma \ref{lem:prime1},
 there is a prime $q$ such that
\begin{enumerate}
\item $q$ is inert in $K$, \item $p$ divides $q+1$ and \item $p$
does not divide  $a_q$.
\end{enumerate}
 And since $q+1\geq p$ by
(2) and $p\geq n+2$, $q$ is strictly greater than $n$. This
implies $q$ does not divide $n$. Therefore, the ring class field
$H_{nq}$ is a proper finite extension of $H_{n}$ and $nq$ is again
a square-free multiple of $r$ and $(nq/r, rN)=1$.

 By the norm-compatibility property given in Proposition
\ref{prop:norm}, when $q\nmid n$ is inert in $K$, we have
$$\mbox{Trace}_{H_{nq}/H_{n}}(P_{nq})=a_qP_{n}.$$

On the other hand,  $E(H_{nq})=E(H_{n})$, by assumption. Hence the
trace of $P_{nq}$ from $H_{nq}$ to $H_{n}$ is divisible by the
degree of $H_{nq}$ over $H_{n}$ which is $q+1$ by (b) of Lemma
\ref{lem:degree}.

Hence, $\mbox{Trace}_{H_{nq}/H_{n}}(P_{nq})$ is divisible by $p$
by the property (3). But by the property (4) and since $p\notin
S$, $p$ divides neither $a_q$ nor the point $P_{n}$, which is a
contradiction. So we have shown that the group generated by the
Heegner points of conductor $rm$ for all square-free integers $m$
is not finitely generated. In particular, this shows that there is
a non-torsion point $P_{nq}\in E(H_{nq})$ but not in $E(H_n)$.

By Lemma \ref{lem:bdd}, for such a non-torsion point $P_{nq}\notin
E(H_n)$, there exists an integer $M$ such that $mP_{nq}\notin
E(H_n)$ for all $m > M$. In other words, the point $P_{nq}$ is
independent of any points in $E(H_n)$. Hence, $E(H_n)\otimes \bbq
\neq E(H_{nq})\otimes \bbq$. Therefore, we conclude that the rank
of $E(H_{rm})$ cannot be bounded, as a square-free integer $m$
such that $(m,N)=1$ goes to infinity.
\end{proof}

Now we prove the unboundedness of the rank of Mordell-Weil groups
over all the ring class fields of conductor $cp^n$ as $n$ goes to
infinity and this will be an important role in proving the main
theorem. To prove this, we need the following simple lemma.

\begin{lem}\label{lem:seq} For an elliptic curve $E/\bbq$ and for a prime $p$ not
dividing the conductor of $E$, let $a_p=p+1-\#E(\bbf_p)$. Then,
there is no infinite sequence $\{c_n\}_{n=0}^{\infty}$ of integers
with $c_0\neq 0$ and satisfying the following linear recurrence,
$$pc_{n+1}=a_pc_n-c_{n-1}, \mbox{ for } n\geq 1,$$
and  for every $N$, there exists $n >N $ such that
 $c_n\neq 0$.
\end{lem}

\begin{proof} Suppose there is such an infinite sequence
$\{c_n\}_{n=0}^{\infty}$ of integers satisfying the above
conditions. Then, the linear recurrence implies that
 $$ (*) \mbox{\hspace{5 cm}} c_n =\alpha^nb_0
 +\beta^nb_1, \mbox{ for
all } n\geq 1,\mbox{\hspace{6 cm}}$$
where $\alpha$ and $\beta$
are two solutions of the quadratic equation
$x^2-\frac{a_p}{p}x+\frac{1}{p}=0$ and
$$b_0=\left(\frac{-\beta}{\alpha-\beta}c_0
\frac{1}{\alpha-\beta}c_1\right) \mbox{ and }
b_1=\left(\frac{\alpha}{\alpha-\beta}c_0
-\frac{1}{\alpha-\beta}c_1\right).$$
 Note that this quadratic
equation has no rational solutions, because if there was, the only
possible pairs of rational solutions are $1$ and $\frac{1}{p}$ or
$-1$ and $-\frac{1}{p}$ and in either case, we get $a_p=\pm (p+1)$
which is impossible since $|a_p|<2\sqrt{p}$ by  Hesse's inequality
(\cite{sil1}, Chapter V, Theorem 1.1). And $b_0\neq 0$ and
$b_1\neq 0$ since $c_0$ is a nonzero integer and $\alpha$ and
$\beta$ are not rational numbers.

Let $F$ be a quadratic extension of $\bbq$ containing $\alpha$ and
$\beta$ and choose an embedding of $F$ into $\overline \bbq_p$
with the valuation $v_p$ such that $v_p(\alpha)<0$ and
$v_p(\beta)=0$, where $\overline\bbq_p$ is an algebraic closure of
the $p$-adic field $\bbq_p$. This is possible because
$\alpha\beta=\frac{1}{p}$ and we can take an automorphism of $F$
and take an embedding of $F$ into $\overline\bbq_p$ such that
$v_p(\alpha)<0$ and $v_p(\beta)=0$.

Since $b_0\neq 0$ and $v_p(\beta)=0$, the recurrence relation
$(*)$ implies that for all large integers $n$,  $v_p(c_n)$ is
dominated by $v_p(\alpha^{n-1})$ which is negative. But for each
$N$, there exists $n>N$ such that $c_n\neq 0$. Since a nonzero
integer $c_n$ has a nonnegative valuation, we get a contradiction.
Hence, there is no such a sequence.
\end{proof}

\begin{prop}\label{prop:inf2} Let $N$ be the conductor of the given
elliptic curve $E/\bbq$. If there is a nontrivial Heegner system
attached to $(E,K)$, where $K$ is a quadratic imaginary  extension
of $\bbq$ which is different from the field $\End(E)\otimes \bbq$,
 then for a prime $p$  such that $p\nmid r\cdot N \cdot [H_r:K]\cdot
disc(H_r)$, there exist an integer $r$ such that the rank
of $E(H_{rp^n})$ is unbounded, as $n$ goes to infinity.
\end{prop}

\begin{proof} Let $p$ be a prime $p\nmid r\cdot N \cdot [H_r:K]\cdot
disc(H_r)$. As we have shown in Proposition \ref{prop:inf}, we
show that  for some integer $r$, the group generated by all of
Heegner points of conductor $rp^n$ for all integer $n\geq 1$ in
the given nontrivial Heegner system is not finitely generated.
Suppose for any integer $r$, this group is finitely generated.
Then, since there are only finitely many torsion (Heegner) points
over all the ring class fields in the system by Lemma
\ref{lem:finite}, for some integer $k$, there is a fixed ring
class field $H_{rp^k}$ over which all Heegner points of conductor
$rp^n$ for all $n \geq 1$ are defined. Let $n_0=rp^k$. Since the
given Heegner system is nontrivial, we may assume that a Heegner
point $P_0$ of conductor $n_0=rp^k$ is of infinite order.

 By the norm-compatibility property given in Proposition
\ref{prop:norm}, 
we have that for all $n\geq 1$,
$$\mbox{Trace}_{H_{n_0p^{n+1}}/H_{n_0p^n}}(P_{n+1})=a_pP_{n}-P_{n-1},$$
where $P_i$ are Heegner points of conductor $n_0p^i$.

Since  $P_{n+1}$ is defined over $H_{n_0}$, hence over
$H_{n_0p^n}$ by assumption, the trace of $P_{n+1}$ from
$H_{n_0p^{n+1}}$ to $H_{n_0p^n}$ is the degree of $H_{n_0p^{n+1}}$
over $H_{n_0p^n}$ which is $p$ by (a) of Lemma \ref{lem:degree}.
Hence, the infinite sequence of Heegner points of conductor of the
form $n_0p^n$ satisfies the linear recurrence relation, $$ (**)
\mbox{ \hspace{4 cm} } pP_{n+1}=a_pP_n-P_{n-1}, \mbox{ for all }
n\geq 1.\mbox{ \hspace{6 cm} }$$

 By the Mordell-Weil Theorem, $E(H_{n_0})$ is finitely generated,
 so by dividing by its torsion subgroup $E(H_{n_0})_{tor}$, all points $P_n \mod
 E(H_{n_0})_{tor}$ lie in $\bbz^k$ for some $k$. Suppose that $E(H_{n_0})\cong
 \bbz Q_1+\cdots +\bbz Q_k+ E(H_{n_0})_{tor}$.
 Now we consider all points $P_n \mod
 E(H_{n_0})_{tor}$ and denote it by $P_n$ again by abuse of
 natation.  Let $P_n=\sum\limits_{i=1}^k c_{n,i}Q_i$ for integers $c_{n,i}$. Since $P_0$ is
 not a torsion point by the assumption, without loss of
 generality we may assume that $c_{0,1}\neq 0$. Then, either
 $c_{2,1}$ or $c_{3,1}$ is nonzero, since otherwise, the relation
 $(**)$ implies that $Q_1$ is a linear combination of points
 $Q_2, Q_3, \ldots,Q_k$ over $\bbz$ which contradicts the linear
 independence of points $Q_i$.  By using $(**)$ and linear dependence of points $Q_i$ inductively, we
 can show that $pc_{n+1}=a_pc_n-c_{n-1}$ for all $n\geq 1$, and
  if $c_{n,1}\neq 0$, then either $c_{n+1,1}\neq 0$ or $c_{n+2,1}\neq
 0$. Hence, by letting $c_n=c_{n,1}$, we get an infinite sequence $\{c_n\}_{n=0}^{\infty}$ of integers
 satisfying the $pc_{n+1}=a_pc_n-c_{n-1}$ for all $n\geq 1$ with
 $c_0\neq 0$ and for each $N$, there exists $n>N$ such that
 $c_n\neq 0$. This is impossible by Lemma \ref{lem:seq}. Therefore, we conclude that the rank of
$E(H_{rp^n})$ is not finitely generated, as $n \rightarrow
\infty$. In particular, by Lemma \ref{lem:bdd}, we conclude
that the rank of $E(H_{rp^n})$ cannot be bounded, as  $n$
goes to infinity.
\end{proof}

Finally, the following proposition proves Theorem \ref{thm:ab}
hence, completes Theorem \ref{thm:main}.

\begin{prop}\label{prop:sigma} Let $\sigma\in \gq$. Let $K$ be a quadratic imaginary
extension of $\bbq$ such that $\sigma|_K\neq id_K$ and $K$ is
different from $\End(E)\otimes \bbq$. Suppose all primes dividing
the conductor $N$ of $E$ are split in $K$. Then, the rank of the
Mordell-Weil group $E(H_n^{\sigma})$ over the fixed subfield of
$H_n$ under $\sigma$ is unbounded, as $n$ goes to $\infty$. Hence,
the rank of $E((K_{ab})^{\sigma})$ is infinite, where $K_{ab}$ is
the maximal abelian extension of $K$.
\end{prop}

\begin{proof} Since all primes
dividing the conductor $N$ of $E$ are split in $K$,  there is a
nontrivial Heegner system attached to $(E,K)$ by (\cite{dar},
Chapter 3, Theorem 3.13). For a given $\sigma\in \gq$, since
$\sigma|_K\neq id_K$, the restriction of $\sigma$ to each ring
class field $H_n$ in the given Heegner system can be lifted as an
involution of $H_n$. Let $\sigma_n=\sigma|_{H_n}$ be the
restriction of $\sigma$ to $H_n$. Then, each ring class field
$H_n$ has a generalized dihedral group structure as its Galois
group over $\bbq$ with an involution $\sigma_n$ such that for any
$\tau\in \Gal(H_n/K)$, $\sigma_n\tau\sigma_n=\tau^{-1}$.

By Proposition \ref{prop:inf2}, we fix a prime $p$ and an integer
$r$ such that $p\nmid r\cdot N \cdot [H_r:K]\cdot disc(H_r)$ and
 the rank of $E(H_{rp^n})$ is
unbounded, as $n$ goes to infinity.  We prove that for an odd
prime $p$ not dividing $rN[H_r:K]disc(H_r)$, the rank of
$E(H_{rp^n}^{\sigma})$ is unbounded as $n$ goes to infinity.
 Suppose not. Then since the restriction $\sigma_n$ of $\sigma$ to $H_{rp^n}$ acts by an involution of each
 ring class field $H_{rp^n}$,
  there exists a fixed integer $n_0=rp^k$ for some $k\geq 1$ such that
$\sigma$ acts by $-1$ on any nontrivial quotient
$(E(H_{n_0p^n})\otimes \bbq)/(E(H_{n_0})\otimes\bbq)$, for all
$n\geq 1$.

By (a) of Lemma \ref{lem:degree}, $\Gal(H_{n_0p^n}/H_{rp}) =
\Gal(H_{rp^{k+n}}/H_{rp})$ is a cyclic group of order $p^{k+n-1}$.
And $\Gal(H_{n_0p^n}/H_{n_0})$ is a subgroup of
$\Gal(H_{n_0p^n}/H_{rp})$. Hence, it is a cyclic subgroup of order
$p^m$ for some $m < k+n-1$.

Let $\tau_n$ be a generator of $\Gal(H_{n_0p^n}/H_{n_0})$.
 Consider $E(H_{n_0p^n})\otimes \bbq$ as
a representation of $\Gal(H_{n_0p^n}/\bbq)$. And for each n, let
$$M_n=(E(H_{n_0p^n})\otimes \bbq)/(E(H_{n_0})\otimes\bbq).$$

For every element $\alpha\in \Gal(H_{n_0p^n}/\bbq)$, let
$\alpha|_{H_{n_0}}$ be the restriction of $\alpha$ to $H_{n_0}$.
Then,  $\alpha|_{H_{n_0}}$ is an element of  $\Gal(H_{n_0}/K)$,
since $H_{n_0}$ is Galois over $\bbq$. Therefore,
$\Gal(H_{n_0p^n}/\bbq)$ acts on $E(H_{n_0})\otimes \bbq$ as well.
So we can consider the quotient $M_n$ as a representation of
$\Gal(H_{n_0p^n}/\bbq)$.

Let $$\rho : \Gal(H_{n_0p^n}/\bbq) \rightarrow \mbox{GL}(M_n)$$ be
the representation of $\Gal(H_{n_0p^n}/\bbq)$. Then, by the
hypothesis, $\sigma_n$ acts by $-1$ on $M_n$. Hence,
$\rho(\sigma_n)=-id$ on $M_n$. On the other hand, by the dihedral
group structure of $\Gal(H_{n_0p^n}/\bbq)$,
$\sigma_n\tau_n\sigma_n=\tau_n^{-1}$. Therefore,
$$\rho(\tau_n^2)=\rho(\tau_n)\rho(\tau_n)=(-id)\rho(\tau_n)(-id)\rho(\tau_n)=\rho(\sigma_n\tau_n\sigma_n\tau_n)=\rho(1)=id.$$

Hence, the restriction of $\rho$ to the cyclic subgroup
$\Gal(H_{n_0p^n}/H_{n_0})$ of $Gal(H_{n_0p^n}/\bbq)$ generated by
$\tau_n^2$ is a trivial representation of $M_n$. Since the order
of $\tau_n$ is an odd integer $p^m$, $$\langle \tau_n^2 \rangle
=\langle \tau_n \rangle = \Gal(H_{n_0p^n}/H_{n_0}).$$ Therefore,
we have
$$M_n^{\Gal(H_{n_0p^n}/H_{n_0})} = M_n^{\langle \tau_n^2 \rangle}=M_n.$$
This implies that
$$E(H_{n_0p^n}^{\Gal(H_{n_0p^n}/H_{n_0})})\otimes\bbq +E(H_{n_0})\otimes\bbq =E(H_{n_0p^n})\otimes\bbq.$$
Since $H_{n_0p^n}^{\Gal(H_{n_0p^n}/H_{n_0})}=H_{n_0}$,
$$E(H_{n_0})\otimes\bbq =E(H_{n_0p^n}) \otimes\bbq, ~~~\mbox{~~for all }
n \geq 1,$$  which is a contradiction to Proposition
\ref{prop:inf2}.

 Hence,
the rank of $E(H_{n_0p^n}^{\sigma})$ is unbounded, as
$n\rightarrow\infty$. Since all ring class fields $H_{n_0p^n}$
are abelian over $K$, this implies that the rank of
$E((K_{ab})^{\sigma})$ is infinite.
\end{proof}


\end{document}